\documentclass[11pt]{article}
\usepackage{ictami}

\newcommand{\PAGENUMBER}{1}     
\title{Quantum Sheaf Cohomology on Surfaces of
General Type I: construction of stable omalous bundles}
\newcommand\shorttitle{Quantum Sheaf Cohomology on Surfaces of
General Type}
\author{C. Anghel }
\begin{document}

\ttl  

\begin{abstract} Quantum sheaf cohomology is a deformation of the cohomology ring of
a sheaf. In recent years, this subject had an impetuous development in
connection with the $(0; 2)$ non-linear sigma model from super-strings theory. The basic piece in this area is a so-called omalous bundle on the variety we start with. After a short overview of the subject, we construct stable omalous bundles  on some classes of surfaces of general type.

\end{abstract}


%
\section{Introduction}

Quantum sheaf cohomology is a deep subject and at the same time a natural extension of the "classical" quantum cohomology developed by Vafa, Witten, Kontsevich and others. However as many subjects that comes from physics, it is defined in general only at a physicist level of rigor. Reecently, Donagi and his collaborators achieved to provide a rigourous foundation of this subject in the case of toric varieties and their tangent bundles \cite{do1}, \cite{do2}. The second Section is devoted to an overview on the subject, following \cite{do1}, \cite{do2} and \cite{guff}. As we shall see, a basic piece in the construction of a quantum sheaf cohomology is a so-called omalous bundle on the variety we start with. The last section is devoted to this subject, namely to the construction of stable omalous bundle on same classes of surfaces. We obtain new examples in rank $2\leq r\leq 6$ on surfaces in $\mathbb P^3$ of high degree. The main technique we use is the Li-Qin construction developed in \cite{liq}. In the sequel of this paper \cite{a4} we will be concerned with computational aspects of the quantum sheaf cohomology for the stable omalous bundles constructed here.

\section{Quantum sheaf cohomology}

\subsection{The "classical" quantum cohomology}

 The quantum cohomology is a formal deformation of the cohomology ring $H^{\ast }(X, \mathbb C)$ for a 
smooth algebraic variety $X$.
 From the physical viewpoint it appears if one try to "observe" the variety $X$ not by its points, but by its rational curves.
 The result is a quantum product on $H^{\ast }(X, \mathbb C)[[\bold q]]$ where ${\bold q}=(q_1,q_2,...)$ is a multi-index whose length is the 
rank of $H^{2 }(X, \mathbb Z)$.
 The main ingredient in the construction is Kontsevich's moduli space of stable maps 
$${\mathcal M}(X,\beta )$$ which parametrize maps ${\mathbb P}^1 \rightarrow X $ whose image has class $ \beta \in H^2(X, \mathbb Z)$. 
 Also, one should mention that the construction depend on a system of marked points, but we shall ignore this aspect here.

 Having the moduli space of stable maps, the construction of the quantum product on $H^{\ast }(X, \mathbb C)$ goes as follows:
 for two classes $\omega_1$, $\omega_2$ $\in H^{\ast }(X, \mathbb C)$ the big deal is to define their quantum product 
$$\omega_1 \star \omega_2 \in H^{\ast }(X, \mathbb C)[[\bold q]].$$
 A first observation is that as consequence of Poincare duality, it is enough to define the pairing 
$$<\omega_1 \star \omega_2, \omega_3 >$$ for any class $\omega_3 \in H^{\ast }(X, \mathbb C)$.
 The pairing can be expanded as a formal sum 
$$\sum_{\beta }<\omega_1 \omega_2 \omega_3 >_{\beta }{\bold q}^\beta $$ 
where $\beta $ varies in the free part of $H^2(X, \mathbb Z) \simeq {\mathbb Z}^r$.
 So, we are reduced to the computation of this triple product $<\omega_1 \omega_2 \omega_3 >_{\beta }$ for each 
$\beta \in H^2(X, \mathbb Z)$.
 This is done by pushing up the three $\omega _i$ in the cohomology of ${\mathcal M}(X,\beta )$ and by taking here the cup product of the resulting classes: 
$$<\omega_1 \omega_2 \omega_3 >_{\beta }=\varphi _{\beta }^{\ast }(\omega_1) \cup  \varphi _{\beta }^{\ast }(\omega_2) \cup  \varphi _{\beta }^{\ast }(\omega_3),$$ where the push up map $\varphi _{\beta } :  H^{\ast }(X, \mathbb C)\rightarrow H^{\ast }({\mathcal M}(X,\beta ), \mathbb C)$ is constructed using the marked points. 
 Even in such a oversimplified picture, we must mention three great difficulties in the full story:
 First of all, the moduli space ${\mathcal M}(X,\beta )$ is not compact at the beginning and its compactification was obtained by Kontsevich. Secondly, this compactification is not a variety but a stack.
 Thirdly, on the compactification one needs a so called virtual fundamental class which is used to give a rigorous meaning for the cup product on ${\mathcal M}(X,\beta )$.

 All these difficulties were resolved for large classes of varieties. Below are two simple examples:
 $$QH^{\ast }(\mathbb P^n)=\frac{\mathbb C[x][[q]]}{(x^{n+1}-q)}$$
 $$QH^{\ast }(\mathbb P^n \times \mathbb P^m)=\frac{\mathbb C[x,y][[p,q]]}{(x^{n+1}-p,y^{m+1}-q)}.$$

\subsection{Quantum cohomology for sheaves}
 A rigorous foundation for the quantum cohomology for sheaves was introduced by Donagi et al. in \cite{do1} and \cite{do2} in connection with the $(0,2)$ nonlinear sigma model.
 Its construction is similar with that for varieties but has a sheaf $E$ on the variety $X$ as supplementary input.
 The sheaf $E$ has to satisfy certain constraint - the omality condition - $$c_1(T_X)=c_1(E) \ \ \ \ c_2(T_X)=c_2(E)$$ which imply the vanishing of the Green-Schwarz anomaly.
 An important point is that the omality is a necessary but not sufficient condition for the existence of a quantum sheaf cohomology for $E$.
 As definition, the quantum sheaf cohomology for $E$ is the structure of a quantum product on $$QH^{\ast }(X,E):=H^{\ast }(X,\Lambda ^{\ast }(E^{\vee}))\otimes \mathbb C[[\bf q]].$$

 It is a deformation of the usual product on the cohomology of $E$.

\begin{remark} For $E=T_X$ the quantum sheaf cohomology of $E$ is the "classical" quantum cohomology of $X$ as one can guess from the  Hodge decomposition. \end{remark}
 The construction goes along the same lines as in the "classical" case: one starts with two elements $\omega_1$, $\omega_2$ $\in H^{\ast }(X, \Lambda ^{\ast }(E^{\vee}))$ and we want to define their quantum product 
$$\omega_1 \star \omega_2 \in H^{\ast }(X, \Lambda ^{\ast }(E^{\vee}))[[\bold q]].$$
 Again by Poinare duality it is enough to define the pairing $$<\omega_1 \star \omega_2, \omega_3 >$$ for any class $\omega_3 \in H^{\ast }(X,\Lambda ^{\ast }(E^{\vee}))$.
 Finally we arrive at the same problem, namely the definition of $<\omega_1 \omega_2 \omega_3 >_{\beta }$ for each 
$\beta \in H^2(X, \mathbb Z)$.
 The next step is to push up the $\omega _i \ 's$ from $H^{\ast }(X, \Lambda ^{\ast }(E^{\vee}))$ to 
$H^{\ast }({\mathcal M}(X,\beta ), \Lambda ^{\ast }(F^{\vee}))$, 
where $F$ is a certain sheaf over ${\mathcal M}(X,\beta )$ obtained from $E$.
 For example, in an ideal situation when the moduli space ${\mathcal M}(X,\beta )$ were fine with classifying map 
$$\varphi : {\mathcal M}(X,\beta )\times \mathbb P^1 \rightarrow X, $$ then $F$ would be $R^0\pi _{1\ast }\varphi {\ast }(E)$.
 Anyway, even in the real world where ${\mathcal M}(X,\beta )$ is not fine, such an $F$ exists and the main point is that
 the omality of $E$ imply $ \Lambda ^{top}(F^{\vee})\simeq K_{{\mathcal M}(X,\beta )}$
 As consequence, if one starts with $\omega _i \in H^{p_i }(X, \Lambda ^{q_i }(E^{\vee}))$ with $\Sigma p_i=dim \ {\mathcal M}(X,\beta )$ and $\Sigma q_i= rank(F)$ then by pushing up the $\omega _i \ 's$ and taking cup product we arrive in $H^{top}({\mathcal M}(X,\beta ),K)\simeq \mathbb C$, producing therefore the desired number $<\omega_1 \omega_2 \omega_3 >_{\beta }$.

\subsection{An example: the quadric surface}

 As one could remark, apart the difficulties connected with the "classical" quantum cohomology, at least to new problems can be seen: the construction of the bundle $F$ and the starting point, namely the construction of the omalous bundle $E$
 In fact, the construction of a quantum sheaf cohomology is known in very few cases. 
 In this section we shall review a result in this direction obtained by Donagi et al. in \cite{do2} concerning the quadric $\mathbb P^1 \times \mathbb P^1$.
 As the starting point on the quadric, the above mentioned authors considers de bundle $E$ as cokernel in the following sequence: $$0\rightarrow {\mathcal O} \oplus {\mathcal O}\rightarrow {\mathcal O}(1,0)^2\oplus {\mathcal O}(0,1)^2\rightarrow E\rightarrow 0,$$
 The left arrow above is given by a matrix of the form
$\begin{bmatrix}
Ax &  Bx \\
Cx' & Dx'
\end{bmatrix}$,
 where $A, B, C, D$ are $2\times 2$ complex matrices and $x=(x_1,x_2)$, $x'=(x'_1,x'_2)$ are homogenous coordinates on the two projective lines.
 Note that the bundle $E$ is a deformation of $T_X$ which correspond to the special values $A=D=I_2$ and $B=C=0$.
 With the previous notations, the result proved by Donagi et al. is: 

\begin{theorem}
For the quadric $X=\mathbb P^1\times \mathbb P^1$ and the sheaf $E$ above, the quantum sheaf cohomology is:
$$QH^{\ast }(X,E)=\frac{\mathbb C[a,b][[p,q]]}{(det(Aa+Bb)=p,det(Ca+Db)=q)}$$.
\end{theorem}
\begin{remark} For $E=T_X$, namely for $A=D=I_2$ and $B=C=0$, the quantum sheaf cohomology above, recover the "classical" quantum cohomology of the quadric, as it is expected. \end{remark}

\subsection{Motivation for the next Section}
 As one can see from this section, the first ingredient for the construction of a quantum sheaf cohomology is an omalous bundle.
 However, despite their importance, only few examples are known in the literature. A first class of examples consists of deformations of the tangent bundle $T_X$ which are omalous by trivial reasons.
 On the other hand, as showed by Andreas and Garcia-Fernandez in 2010, the stability of an omalous bundle is a very important property, because such a bundle provide a solution of the so-called Strominger system in super-string theory.
 In this direction, a first systematic attempt to construct stable omalous bundles was done in 2011 by Henni and Jardim in \cite{hj}. They uses monads to construct:\\
- stable omalous rank $3$ bundles on $3$-folds in $\mathbb P^4$,\\
- stable omalous rank $2$ bundles on c.i. CY's in projective spaces,\\
- omalous of rank $>3$ bundles on multi-blowup of the plane,\\
- stable omalous of various ranks on the Segre variety.
 Also, in \cite{am}, was  studied omalous rank $2$-bundles on Hirzebruch surfaces.
 The next section will be devoted to the construction of stable omalous bundles on surfaces, with special emphasis for the case where $X$ is a surface of general type.
 I shall first describe a general construction for stable bundles due to Li and Qin \cite{liq}.
After that I shall consider certain concrete cases where this construction can be applied to produces stable omalous bundles.

\section{Stable omalous bundles on surfaces of general type}
\subsection{General construction of stable bundles on surfaces}

In what follows $X$ is a smooth projective surface and $L$ a very ample polarization. The $L$-stability of a sheaf $E$ means that it has the greatest fraction $$\frac{c_1\cdot L}{rank}$$ among all its sub-sheaves.
 The main problem for the moment is the following: fixing the rank $r$ and the first Chern class $c_1$, to find a computable bound $\alpha $ depending only on $r$, $L$ and $c_1$ such that for any $c_2\geq \alpha $ there is an $L$-stable vector bundle of rank 
$r$ with the given Chern classes $c_1$ and $c_2$. 
 The main result of Li and Qin asserts that for $\alpha $ one can take the following value:
 $$\alpha = (r-1)[1+ max(p_g,h^0(S,{\mathcal{O}}_X(rL-c_1+K_X)))+4(r-1)^2\cdotp L^2]+$$ 
$$+(r-1)c_1\cdotp L-\frac{r(r-1)}{2}\cdotp L^2,$$ 
 where $K_X$ is the canonical class and $p_g=h^0(X,K_X)$ the geometric genus of $X$.
 The main point in the proof of Li-Qin theorem is the following generalization of the usual Cayley-Bacharach property cf. \cite{liq}:
\begin{lemma}
Consider $r-1$ line bundles $L_1,...,L_{r-1}$ and $0$-cycles \\ $Z_1,...,Z_{r-1}$ on $X$; 
let $W=\oplus( {\mathcal O}_X(L_i)\otimes {\mathcal I}_{Z_i})$. Then, there is a locally free extension in $Ext^1(W,{\mathcal O}_X(L'))$ iff for any $i=1...(r-1)$, $Z_i$ satisfies the Cayley-Bacharach property with respect to the linear system ${\mathcal O}_X(L_i-L'+K_S)$. 
\end{lemma}
After that, the desired bundle $E$ is constructed as an extension 
$$0\rightarrow {\mathcal O}_X(c_1+(1-r)L)\rightarrow E \rightarrow \oplus( {\mathcal O}_X(L)\otimes {\mathcal I}_{Z_i})\rightarrow 0,$$
for a convenient choice of $r-1$ reduced  $0$-cycles $Z_i's$ that ensures the stability of $E$. A similar construction in rank $2$ appeared in \cite{A}, \cite{A1} and was used in \cite{a3} to construct large families of bundles in connection with the strong Bogomolov inequality. 

 As a conclusion with respect to the Li-Qin construction we can state the following:
\begin{remark} The value of $\alpha $ grows up with the Chern number $c^2_1$ due to the presence of the $h^0$-term and to Riemann-Roch.
\end{remark}
 Therefore, if their construction can produce omalous bundles, it is better to try on surfaces that satisfy at least $$c_2>>c^2_1.$$
 The above inequality, combined with Bogomolov-Miyaoka-Yau inequality $c^2_1\leq 3c_2$ for surfaces of general type, suggests to look at certain convenient such surfaces.

\subsection{Stable omalous bundles on "good" surfaces of general type}

 Viewing the above considerations, we can introduce the following:
\begin{definition}
A surface of general type is "good of type $(r,L)$" if for $c_1=\pm K_X$ and $c_2=c_2(X)$, there exists $r\in \mathbb N$ and a very ample line bundle $L$, such that $$c_2\geq \alpha (r,c_1,L),$$ where $\alpha $ is the Li-Qin constant introduced before.
\end{definition}
 In terms of the above definition, the Li-Qin existence Theorem has the following obvious consequence:
\begin{corollary}
On a general type surface $X$ "good of type $(r,L)$", there exists $L$-stable omalous vector bundles of rank $r$.
\end{corollary}
\subsection{Examples of "good" surfaces}

This last part is devoted to the illustration of the above results on a concrete class of examples: $X_d$ a smooth surface of degree $d$ in $\mathbb P^3$.
 Well known computations gives for $X_d$ the following values of invariants: 
$$c_2=d^3-4d^2+6d \ \ \ \ c^2_1=d(d-4)^2$$
 Also, by Noether formula, $$p_g=\chi ({\mathcal O})-1=\frac{c^2_1+c_2}{12}=\frac{d^3}{6}+....$$
 So, the leading therm in $d$ which appear in the Li-Qin constant $\alpha $ is $$\frac{d^3(r-1)}{6}$$
 As consequence, for $ 2 \leq r \leq 6 $ and $ d>>0 $ we have $ c_2 \geq \alpha $. 
 So, on can apply the Li-Qin construction, obtaining the following:
\begin{theorem}
There is an explicitly computable constant $d_0$, such that for all $ d \geq d_0 $ and all $ 2 \leq r \leq 6 $, on any smooth surface 
$ X_d \subset \mathbb P^3 $ of degree $d$ there exists a stable omalous bundle of rank $r$.
\end{theorem}

\section{Further directions}
 As further directions, an open question asked by Donagi et al. in \cite{do3} concern the computation of the quantum sheaf cohomology for other sheafs than the tangent bundle.
 Of course, the above constructed stable omalous bundles are good candidates for such a computation.
 Moreover, due to the range of their ranks, less than $7$, this question is very interesting viewing the following conjecture from \cite{guff}:

"For omalous bundles $E$ of rank $r\leq 7$ on a smooth variety, the quantum sheaf cohomology $QH^{\ast }(X,E)$ exists".

%
%
\auth                                                     
     Cristian Anghel\\
     Department of Mathematics\\
     Institute of Mathematics of the Romanian Academy\\
     Bucharest, Romania \\
     email: {\em Cristian.Anghel@imar.ro}

%
%
\end{document}